\documentclass[12pt]{article}
\usepackage{graphicx, epsfig,amsmath,amssymb}
\usepackage{e-jc}
\specs{R14}{16}{2009}

\numberwithin{equation}{section} 

\newcommand{\abs}[1]{\left\vert#1\right\vert}

\begin{document}

\date{\dateline{Feb 2, 2008}{Jan 17, 2009}{Jan 23, 2009}}

\title{On the monochromatic Schur Triples type problem}

\author{Thotsaporn ``Aek'' Thanatipanonda\\
\small{Department of Mathematics and Computer science}\\[-0.8ex]
\small{Dickinson College}\\[-0.8ex]
\small{Carlisle, PA 17013, USA}\\
\small{\tt thanatit@dickinson.edu}%
}

\maketitle

\begin{abstract}

We discuss a problem posed by Ronald Graham about the minimum number,
over all 2-colorings of $[1,n]$, of monochromatic $\{x,y,x+ay\}$ triples for
$a \geq 1$. We give a new proof of the original case of $a=1$.
We show that the minimum number of such triples is at most
$\frac{n^2}{2a(a^2+2a+3)} + O(n)$ when $a \geq 2$. We also find a new upper bound for the
minimum number, over all $r$-colorings of $[1,n]$, of monochromatic Schur
triples, for $r \geq 3$.

\end{abstract} %

\section{Introduction}

The {\it Schur numbers}, $s(r)$, denote the maximal integer $n$
such that there exists an $r$-coloring of $[1,n-1]$ that avoids a
monochromatic solution to $x+y=z$. For example $s(2) = 5$ and $s(3) = 14$.
$s(5)$ is unknown but is conjectured to be 161.

The original question
about the minimum number, over all 2-colorings of $[1,n]$, of
monochromatic Schur triples was asked by Ronald Graham in 1997. It can be
thought of as a bigger scale version of Schur numbers. It was solved
in 1998. The answer is $\frac{n^2}{22}+O(n)$
that is realized by coloring  the first
$\frac{4n}{11}$ integers red, the next $\frac{6n}{11}$ integers blue,
and  the final $\frac{n}{11}$ integers red.
The first two solutions were given by Robertson and  Zeilberger \cite{RZ}
and Schoen \cite{Schoen}. Later Datskovsky \cite{BD} found another proof.

Ronald Graham asked another question generalizing the
original one. The question was about the minimum number of
monochromatic
$(x,y,x+ay)$ triples, $a \geq 2$ on $[1,n]$.
We discuss this problem in this paper.

In Section 2, we give a new simple proof of
the original problem of finding the minimum number,
over all 2-colorings of $[1,n]$,
of monochromatic Schur triples. In Section 3,
we talk about the generalized problem asked by Graham.
For this problem,
we wrote a computer program to find an optimal
coloring for small $n$ to see some patterns. Then
we used a newly found ``greedy calculus" to
obtain a ``good" upper bound.
The final step was to try to match the lower bound
and upper bound of the problem.
In Section 4, we also apply the greedy
calculus to the original question
on Schur triples with $ r \geq 3$, to obtain
a new upper bound.

\section{The minimum number, over all 2-colorings of $[1,n]$, of monochromatic Schur triples}

\noindent \textbf{2.1 A Greedy Algorithm for The Upper bound} \\

It is natural to find examples of good colorings first.
This example will give us an upper bound. Then we try to show that this upper bound is also a lower bound.

We will show how to find an upper bound for the minimum number, over all
2-colorings of $[1,n]$, of monochromatic triples that are solutions of
$x+y=z$. We will obtain this upper bound
by using the Greedy Algorithm. We denote the colors red and blue.

The general idea is to keep adding more new intervals with
different colors so that, each time, the overall
coloring has the least number of monochromatic triples.
For other proofs of this original problem see \cite{RZ}, \cite{BD}, \cite{Schoen}. \\

\noindent \underline{\textbf{First}}\\

\noindent We paint the first interval of
length $k$ red. We will have $\frac{k^2}{4}$ monochromatic triple
solutions of $x+y=z$ (we are assuming $x \leq y$). \\

\noindent \textbf{Note:} $O(k)$ terms are suppressed in this exposition. \\

\noindent \underline{\textbf{Second}}\\

\noindent We paint the second interval blue. We want to find the length of the
interval (with this color) so that the overall number of the monochromatic
triples is minimized. \\

\noindent Let the length of this interval be
$(1+j)k$ (here $j$ is the number we want to find).\\

\noindent The total number of monochromatic triples on the whole interval is now
$\frac{k^2}{4}+ \frac{j^{2}k^{2}}{4} = \frac{(1+j^2)k^2}{4}$.
The total length is n = $k+(1+j)k = (2+j)k.$ \\ \noindent So
the total number of monochromatic triples in terms of $n$ is
$\frac{(1+j^2)(\frac{n}{2+j})^2}{4} =
\frac{(1+j^2)}{(2+j)^2}\frac{n^2}{4}$.\\

\noindent To find the minimum,
we use calculus to get $j = \frac{1}{2}$. The total number of
monochromatic Schur triples is then
$\frac{n^2}{20}+O(n)$.\\

\noindent So far so good. We have a coloring that paints the first $k$
integers red, followed by painting the next $(1+\frac{1}{2})k$
integers blue.\\

\noindent \underline{\textbf{Third}}\\

\noindent Now we try to stick red at the end of the interval,
and try to
lower the overall number of triples. Say the length of this interval is
$jk$, where $j$ is the number we want to find.The total
length is $n =k + (1+\frac{1}{2})k + jk = (\frac{5}{2}+j)k.$ \\

\noindent \textbf{Case 1:} $j\leq 1$ \\
\noindent The total number of monochromatic Schur triples on the whole interval is
$\frac{k^2}{4}+ \frac{k^2}{16}+\frac{j^{2}k^{2}}{2} =
(\frac{5}{16}+\frac{j^2}{2})k^2$. \\
\noindent So the total number of monochromatic Schur triples in terms of $n$
is\\
$(\frac{5}{16}+\frac{j^2}{2})\frac{n^2}{(\frac{5}{2}+j)^2} = \frac{5+8j^2}{(5+2j)^2}\frac{n^2}{4} $.\\
To find the minimum, we again use calculus and get $j = \frac{1}{4}$. The
total number of monochromatic triples in this case is $\frac{n^2}{22}+O(n)$. \\

\noindent\textbf{Case 2:} $ 1 \leq j \leq \frac{5}{2}$ \\
\noindent The total number of monochromatic Schur triples on the whole interval is
$\frac{k^2}{4}+ \frac{k^2}{16}+(j-\frac{1}{2} )k^2=
(j-\frac{3}{16})k^2$. \\
\noindent So the total number of monochromatic Schur triples in terms of $n$ is
$(j-\frac{3}{16})\frac{n^2}{(\frac{5}{2}+j)^2}$.\\
We again use calculus to find the minimum. We get $j = 1$. The
total number of monochromatic triples in this case is $\frac{13}{196}n^2+O(n)$. \\

\noindent \textbf{Case 3:} $ \frac{5}{2} \leq j $ \\
\noindent The total number of monochromatic Schur triples on the whole interval is
$\frac{k^2}{4}+ \frac{k^2}{16}+ (j-\frac{1}{2} )k^2 + \frac{(j-\frac{5}{2})^2k^2}{4}$. \\
\noindent The total number of monochromatic Schur triples in terms of $n$ is
$(2j^2-2j+11)\frac{n^2}{8(\frac{5}{2}+j)^2}$.\\
We again use calculus to find the minimum. We get $j = \frac{5}{2}$. The
total number of monochromatic triples in this case is $\frac{37}{400}n^2+O(n)$. \\

\noindent In conclusion, the total minimum is $\frac{n^2}{22}+O(n)$.
The coloring for the whole interval is a red interval of length equal
to $k$, a blue interval of length equal to $(1+\frac{1}{2})k$ and
another red interval of length equal to $\frac{1}{4}k$. $k$ is such that the sum of these intervals is $n$, i.e.
$k=\frac{n}{(\frac{5}{2}+\frac{1}{4})}= \frac{4n}{11}$. \\

\noindent \underline{\textbf{Fourth}}\\

\noindent We try to lower the
bound further by having a blue interval of length, say, $jk$ at the end of the
previous interval. But now we get that the minimizing
$j$ is negative. So we stop.\\

\noindent As a conclusion, the optimal coloring with respect to the
greedy algorithm is proportional to
$[1,\frac{3}{2},\frac{1}{4}]$, with colors $[R,B,R]$
yielding that indeed the minimal number is $\frac{n^2}{22}+O(n).$
 \\ %

\noindent \textbf{2.2 The Lower Bound} \\

Finding a lower bound is, in general, the difficult part. However, in this case, it is possible
 since we can turn the problem into a calculus problem. A similar technique was used in \cite{RPS}.\\

\noindent \textbf{Definition}\\
Let $M_\chi(n)$ be the number of monochromatic Schur triples for a 2-coloring $\chi$ of $[1,n]$.\\
Let $Q$ be twice the number of non-monochromatic Schur triples for a $2$-coloring of $[1,n].$\\

\noindent Divide the interval $[1,n]$ into $k$ consecutive intervals. \\
Let $r_i$ be the number of red points in the interval $I_i$.\\
Let $b_i$ be the number of blue points in the interval $I_i$.\\
Let $S_{i,j}$ be the number of non-monochromatic pairs in the square of $I_i \times I_j.$\\
Let $T_{i,j}$ be the number of non-monochromatic pairs
in the triangle of $I_i \times I_j$.

\medskip

\noindent \textbf{Note:}$ \;\ r_i+b_i = \frac{n}{k}$. \\

\noindent \textbf{Lemma 1)}
$M_{\chi}(n)=\frac{n^2}{4}-\frac{Q}{2}+O(n).$ \\

\noindent {\it Proof: \ }
The total number of triples is
$$\abs{ \mbox{monochromatic triples}} + \abs{ \mbox{non-monochromatic triples}}
= M_{\chi}(n) + \textstyle\frac12{Q}.$$
Since the total number of triples is
$\frac{n^2}{4}+O(n)$, we have $ M_{\chi}(n) =
\frac{n^2}{4}-\frac{Q}{2}+O(n). \;\ \Box$ \\

\noindent The plan is to find an upper bound of $Q$ that will give a lower bound for $M_{\chi}(n)$. \\

\noindent \textbf{Lemma 2)} $\displaystyle Q = \abs{R}\abs{B} +\frac{1}{2}(\sum_{i +j < k}S_{i,j} + \sum_{i=1}^kT_{i,k-i+1}), $
where $\abs{R} = \displaystyle\sum_{i=1}^kr_i$ and
$\abs{B} = \displaystyle\sum_{i=1}^k b_i$.\\

\noindent {\it Proof: \ }
\[\begin{array}{llll}
Q \;\ &=& \;\ \abs{\{(R,B),(B,R)| \;\ y-x \geq 0\}} +
\abs{\{(R,B),(B,R)| \;\ x+y \leq n, x \geq y \}} \\
&=& \;\ \abs{\{(R,B),(B,R)| \;\ y-x \geq 0\}} + \frac{1}{2}\abs{\{(R,B),(B,R)| \;\ x+y \leq n \}}. \\
\end{array} \]

\noindent Note that each non-monochromatic triple contributes two non-monochromatic pairs: \;\
for example, $(x,y,z)=(R,B,R)$ gives $(x,y)=(R,B)$ and $(y,z)=(B,R)$. The statement of the lemma follows. $\Box$\\

\noindent Now we find an upper bound for $Q$.
For each $T_{i,j}$ we have two ways to bound it: \\
1) $T_{i,j} \leq$ area of the triangle = $\frac{1}{2}(\frac{n}{k})^2.$ \\
2) $T_{i,j} \leq S_{i,j}$. \\

\noindent \textbf{Example 1:} $k = 2$, with the upper bound of $T_{1,2}, T_{2,1}$
using the areas of the triangles.
We have
\[\begin{array}{llll}
Q  &=&  \abs{R}\abs{B} + \frac{1}{2}(S_{1,1} + T_{1,2} + T_{2,1}).  \\
&\leq& (r_1+r_2)(b_1+b_2) + r_1b_1+\frac{n^2}{8}. \\
 &=& (r_1+r_2)(n-r_1-r_2) + r_1(\frac{n}{2}-r_1)+\frac{n^2}{8}. \\
\end{array} \]

\noindent We use calculus to find a maximum of $Q$ where $0 \leq r_1, r_2 \leq \frac{n}{2}$.
The optimal solutions is $r_1 = \frac{n}{4}$ and $r_2= \frac{n}{4}$.\\
We then get the maximum $Q$ as $\frac{7n^2}{16}$. This yields $M_\chi(n) \geq \frac{n^2}{32}+O(n).\;\ \Box$  \\

\noindent \textbf{Example 2:} $k = 3$, with the upper bound of
$T_{1,3}, T_{3,1}$ using the areas of the triangles and
the upper bound of $T_{2,2}$ using $S_{2,2}$. We have

\[\begin{array}{llll}
Q &=&  \abs{R}\abs{B} + \frac{1}{2}(S_{1,1} + S_{1,2}+ S_{2,1}+ T_{1,3} + T_{2,2}+ T_{3,1}).  \\
&\leq& (r_1+r_2+r_3)(b_1+b_2+b_3) + r_1b_1 + r_1b_2 + r_2b_1+  r_2b_2+\frac{n^2}{18}. \\
\end{array} \]

\noindent We use calculus to find a maximum of $Q$ where $0 \leq r_1, r_2, r_3 \leq \frac{n}{3}$.
One of the optimal solution is $r_1 = 0, r_2 = \frac{n}{3}$ and $r_3= \frac{n}{6}$.\\
This yields the maximum $Q$ is $\frac{5n^2}{12}$ which leads to $M_\chi(n) \geq \frac{n^2}{24}+O(n).\;\ \Box$  \\

\noindent This is pretty nice. We can use calculus to get a decent
lower bound of the problem. The calculation can even be done by hand.
The hope to match the upper bound and lower bound is to try 11 intervals.
This time we need a computer to help doing the calculation.\\

\noindent \textbf{Example 3:} $k=11$,  \\
We bound $T_{2,10}, T_{3,9},T_{4,8}, T_{8,4}, T_{9,3}$ and $T_{10,2}$ by the area of each triangle which is $\frac{n^2}{242}$.\\
We bound $T_{i,12-i}$ by $S_{i,12-i}$, where $i =1,5,6,7,11$.\\

\noindent We get eight optimal solutions to the maximum of $Q$. One of them is \\
$[r_1, r_2,\cdots, r_{11}]
= [\frac{n}{11},\frac{n}{11},\frac{n}{11},\frac{n}{11},0,0,0,0,0,0,\frac{n}{11}].$ \\
This yields the maximum of $Q$ as $\frac{9n^2}{22}$ which gives $M_\chi(n) \geq \frac{n^2}{22}+O(n). \;\ \Box$  \\

\noindent Since the lower bound matches the upper bound, the problem is solved.


\section{Generalized problem, $x+ay=z$, $a \geq 2$}

\noindent \textbf{3.1 A Greedy Algorithm for Upper bounds} \\

We will show how to find an upper bound for the minimum number, over all
2-colorings of $[1,n]$, of monochromatic triples that are solutions of
$x+ay=z$, for a fixed integer $a \geq 2$ (we are no longer stipulate $x \leq y$).
We will obtain this upper bound
by using the Greedy Algorithm. The general idea is the same as in the previous section.
We again call the colors red and blue.\\

\noindent \underline{\textbf{First}}

\medskip

\noindent We paint the first interval of
length $k$ red. We will have $\frac{k^2}{2a}$ monochromatic triples as
solutions of $x+ay=z$.\\

\noindent \underline{\textbf{Second}}

\medskip

\noindent We paint the second interval blue. We want to find the length of the
interval (with this color) so that the overall number of monochromatic
triples is minimum. \\

\noindent Let the length of this interval be
$(a+j)k$ (here $j$ is the number we want to find).\\

\noindent The total number of monochromatic triples on the whole interval is now
$\frac{k^2}{2a}+ \frac{j^{2}k^{2}}{2a} = \frac{(1+j^2)k^2}{2a}$. \\
The total length $n$ is $k+(a+j)k = (1+a+j)k.$ \\
\noindent So the total number of monochromatic triples in terms of $n$
is $\frac{(1+j^2)(\frac{n}{1+a+j})^2}{2a} =
\frac{(1+j^2)}{(1+a+j)^2}\frac{n^2}{2a}$.\\

\noindent To find the minimum,
we use calculus to get $j = \frac{1}{a+1}$. The total number of
monochromatic triples is then
$\frac{n^2}{2a(a^2+2a+2)}$.\\

\noindent So far so
good. We have a coloring  that paints  the first $k$ integers red, followed
by painting the next $(a+\frac{1}{a+1})k$ integers blue.\\

\noindent \underline{\textbf{Third}}

\medskip

\noindent Now we try to stick red at the end of the interval,
and try to
lower the overall number of triples. Say the length of this interval is
$jk$, where $j$ is the number we want to find. The total
length $n$ is $k + (a+\frac{1}{a+1})k + jk = (1+a+\frac{1}{a+1}+j)k.$ \\

\noindent \textbf{Case 1:} $j \leq a$ \\
The total number of monochromatic Schur triples on the whole interval is
$\frac{k^2}{2a}+ \frac{k^2}{2a(a+1)^2}+\frac{j^{2}k^{2}}{2a}.$ \\
\noindent So the total number of monochromatic Schur triples in terms of $n$
is\\ $
(\frac{1}{2a}+ \frac{1}{2a(a+1)^2}+\frac{j^{2}}{2a})(\frac{n}{(1+a+\frac{1}{a+1}+j)})^2$.\\

\noindent To find the minimum, we again use calculus to get $j = \frac{1}{a+1}$. The
total number of monochromatic triples is $\frac{n^2}{2a(a^2+2a+3)}$.  \\

\noindent \textbf{Case 2:} $ j \geq a $ \\
The total number of monochromatic Schur triples on the whole interval is
at least
$\frac{k^2}{2a}+ \frac{k^2}{2a(a+1)^2}+\frac{k^{2}}{2a}+(j-a)k^2.$ \\
\noindent The total number of monochromatic Schur triples in terms of $n$
is

$(\frac{1}{2a}+ \frac{1}{2a(a+1)^2}+\frac{1}{2a}+(j-a))(\frac{n}{1+a+\frac{1}{a+1}+j})^2.$

\noindent To find the minimum, we again use calculus to get $j = \frac{3a^4+7a^3+4a^2-2a-3}{a(a+1)^2}$.
The total number of monochromatic triples is $\frac{a(a+1)^2n^2}{2(4a^4+10a^3+8a^2-3)}$.  \\

\noindent The total number of triples in case 2 is always bigger than
the one in case 1 for $a \geq 2$. In conclusion, the minimum total number of monochromatic
triples relative to this method is $\frac{n^2}{2a(a^2+2a+3)}$.
The coloring for the whole interval is a red interval of length equal
to $k$, a blue interval of length equal to $(a+\frac{1}{a+1})k$ and
another red interval of length equal to $\frac{1}{a+1}k$. $k$ is such that
the sum of these intervals is $n$, i.e.
$k=\frac{n}{(1+a+\frac{2}{a+1})}$. \\

\noindent \underline{\textbf{Fourth}}

\medskip

\noindent We try to lower the
bound even further by having a blue interval of length, say, $jk$ at the end of the
previous interval. But now we get that the minimizing
$j$ is negative. So we stop.\\

\noindent As a conclusion, the optimal coloring is proportional to
$[1,a+\frac{1}{a+1},\frac{1}{a+1}]$, with colors $[R,B,R]$
yielding that indeed the minimal number is $\frac{n^2}{2a(a^2+2a+3)}+O(n).$
 \\ %

\noindent \textbf{3.2 Lower bounds} \\

\noindent We will use a similar technique for the lower bound of the original problem.
We find an upper bound for non-monochromatic triples in $[1,n]$.
This gives a lower bound for the number of monochromatic triples. \\

\noindent We use the notation $(R,B)$ and $(B,R)$ for the
non-monochromatic pair $(x,y)$.\\

\noindent \textbf{Definition:} \\
Let $\abs{R}$ be the number of red points in $[1,n]$.\\
Let $\abs{B}$ be the number of blue points in $[1,n]$.\\

\noindent \textbf{Lemma 3)}
$\abs{\{(R,B),(B,R) | \;\ y > x, \;\ y-x \;\ \mbox{is divisible by } a \}}
\leq \frac{\abs{R}\abs{B}}{a}$. \\

\noindent {\it Proof: \ }
Let $\abs{r_{i}}$ = number of red points in the congruence class $i$ (mod $a$).  \\
Let $\abs{b_{i}}$ = number of blue points in the congruence class $i$ (mod $a$).

\noindent We remark that
$r_{i} + b_{i} = \frac{n}{a}$, $1 \leq i \leq a$ and
$\sum_{i=1}^a r_{i} = \abs{R}$.
\begin{align*}
\big|\{(R,B),(B,R) |& \;\ y > x, y-x \;\ \mbox{is divisible by }a \}\big|
- \frac{\abs{R}\abs{B}}{a} \\
& = \sum_{i=1}^a r_{i}b_{i} - \frac{\abs{R}\abs{B}}{a}  \\
& = \sum_{i=1}^a r_{i}(\frac{n}{a}-r_{i})
    -\frac 1a\Big(\sum_{i=1}^a r_{i}\Big)\Big(n-\sum_{i=1}^a r_{i}\Big)  \\
& =  -\sum_{i=1}^a r_{i}^2 + \frac 1a(\sum_{i=1}^a r_{i})^2  \\
& \leq 0,\mbox{ \ by the Cauchy-Schwarz inequality.}
\end{align*}
Moreover, equality holds when $r_{1} = r_{2} =\cdots = r_{a}$. \;\ $\Box$ \\

\noindent Let $Q_a$ be two times the number of
non-monochromatic triples of solutions of
$x+ay=z$ in a $2$-coloring of $[1,n]$.\\

\noindent \textbf{Lemma 4)} \;\ $Q_a \leq
\frac{\abs{R}\abs{B}}{a}
+ \big|\{(R,B),(B,R)| \,\ y-ax \geq 0\}\big|
+ \big|\{(R,B),(B,R)| \,\ y+ax \leq n \}\big|\!.$

\noindent {\it Proof: \ }
\[\begin{array}{llll}
Q_a &=& \big| \{ \mbox{the non-monochromatic pair} \;\ (x,y) |
\ y > x \;\ \mbox{and} \;\ y-x \;\ \mbox{is divisible by a} \}\big| \\[2pt]
& & + \big|\{ \mbox{the
non-monochromatic pair} \;\ (x,y) | \;\ y-ax \geq 0\} \big| \\[2pt]
& & + \big|\{
\mbox{the non-monochromatic pair} \;\ (x,y) | \;\ y+ax \leq n \}\big| \\[1ex]
& \leq & \frac{\abs{R}\abs{B}}{a} + \abs{\{(R,B),(B,R) | \;\ y-ax \geq 0\}} +
\big|\{(R,B),(B,R) | \;\ y+ax \leq n \}\big|\\ & & \mbox{by Lemma 3.} \;\ \Box
\end{array} \] \\

\noindent When the points on the
$x$-axis and the $y$-axis are painted
with either color red or blue,
$\abs{\{(R,B),(B,R) | \;\ y+ax \leq n\}}$ is
the number of non-monochromatic coordinate pairs inside the triangle 1 below. \\
Similarly $\abs{\{(R,B),(B,R) | \;\ y-ax \geq 0\}}$ is the number of
non-monochromatic coordinate pairs inside the triangle 2. \\

\begin{figure}[htb]
\[\includegraphics{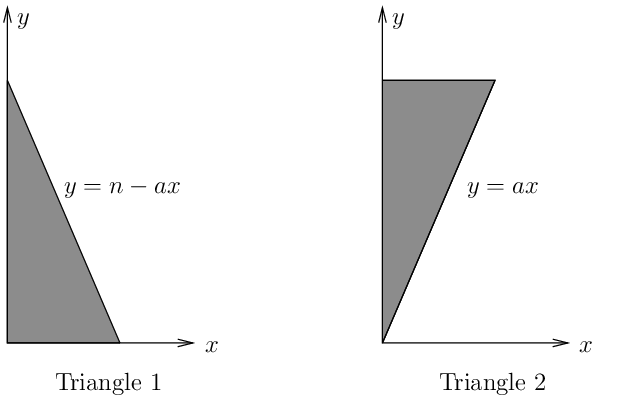}\]
\end{figure}

\noindent Divide the interval $[1,n]$ into $k$ consecutive intervals. \\
Let $r_i$ be the number of red points in the interval $I_i$.\\
Let $b_i$ be the number of blue points in the interval $I_i$.\\
Let $S_{i,j}$ be the number of non-monochromatic pairs in the square $I_i \times I_j.$\\
Let $T_{i,j}$ be the number of non-monochromatic pairs
in the intersection of each of the triangle we consider and the square $I_i \times I_j$.\\

\noindent \textbf{Note:}$ \;\ r_i+b_i = \frac{n}{k}$. \\

\noindent \textbf{Theorem 1)} $Q_2 \leq \frac{57n^2}{121}+O(n)$.\\

\noindent We find an upper bound on $Q_2$
by using calculus on the equation from the previous lemma. \\

The main part of calculating $Q_2$ is to compute the maximum number of
non-mono\-chromatic pairs in triangle 1 and triangle 2 in the pictures
above.  However there are $I_i \times I_j$ for some $i,j$ that
intersect the triangle only partly. We denote them $T_{i,j}$. \\

\noindent For each $T_{i,j}$, in the triangle we consider, we have two ways to bound it,  \\
1) $T_{i,j} \leq$ area of the intersection of triangle and the square $I_i \times I_j$.\\
2) $T_{i,j} \leq S_{i,j} = r_ib_j+r_jb_i$. \\

\noindent In this case, we use 11 intervals, $k=11$.  \\

\noindent In triangle 1, we bound $T_{1,10}$, $T_{2,9}$, $T_{2,8}$, $T_{3,7}$,
$T_{3,6}$, $T_{5,3}$, $T_{5,2}$ and $T_{6,1}$
by the area of each intersecting triangle. We bound $T_{1,11},T_{4,5}$ and $T_{4,4}$ by $S_{i,j}$.\\
In triangle 2, we bound $T_{2,4}, T_{3,5},T_{3,6}$ and $T_{6,11}$ by the area of each intersecting triangle.
We bound $T_{1,1},T_{1,2},T_{2,3},T_{4,7},T_{4,8},T_{5,9}$ and $T_{5,10}$ by $S_{i,j}$.\\
We then run the Maple program.
We get four optimal solutions to the maximum of $Q_2$. Two of them are
$[r_1, r_2,\dots, r_{11}] = [\frac{n}{11},\frac{n}{11},0,\frac{n}{11},0,0,0,0,0,\frac{n}{11},0]$
and  \\
$[\frac{n}{11},\frac{n}{11},0,0,0,0,0,0,\frac{n}{11},\frac{n}{11},\frac{n}{11}]$.
The other two are the switching colors of the first two. \\
This yields an upper bound on $Q_2$ of $\frac{57n^2}{121}+O(n). \;\ \Box$  \\

\noindent \textbf{Definition:}

\smallskip

\noindent Let $M_{\chi,a}(n)$ be the number of monochromatic triples of solutions of $x+ay=z$
for a 2-coloring $\chi$ of $[1,n]$.\\

\noindent \textbf{Corollary} $M_{\chi,2}(n) \geq \frac{7n^2}{484} +O(n).$ \\

\noindent {\it Proof: \ }
The total number of triples is
$$\abs{ \mbox{monochromatic triples}} + \abs{ \mbox{non-monochromatic triples}}
= M_{\chi,a}(n) + \frac{Q_a}{2}.$$
Since the total number of triples is
$\frac{n^2}{2a}+O(n)$, we have $ M_{\chi, a}(n) \geq
\frac{n^2}{2a}-\frac{Q_a}{2}+O(n)$.
The lower bound on $M_{\chi,2}(n)$ follows from the upper bound on
$Q_2$ from Theorem 1. $\Box$ \\

\noindent \textbf{Note:} \\

\noindent 1) For $a = 3$, we found, $M_{\chi,3}(n) \geq \frac{n^2}{2268}+O(n)$. We ran the calculus program on 9 intervals
with a particular upper bound of $T_{i,j}$.\\

\noindent 2) For case $a \geq 4$, we could not find a positive lower bound for $M_{\chi,a}(n)$ yet.
One of the reasons is that the upper bound of $M_{\chi,a}(n)$ is very small.\\

\section{The minimum number, over all $r$-coloring of $[1,n]$, of monochromatic Schur triples}

\noindent \textbf{4.1 A Greedy Algorithm for The Upper bounds} \\

The method to obtain the
upper bounds in this section is similar to the one used in sections 2 and 3. In
general we start with the first interval having color 1. Then we add
interval 2 with color 2 in the optimal way. Then we add the third interval
starting with color 1. If we get a positive solution, we move to the
fourth interval. Otherwise we try with color 3. We keep going on in this
fashion until there is no color that gives a positive solution.

Since there are many intervals involved in the computation, it is too much
computation to do by hand. We wrote a computer program to help us compute
the solutions for each $r$-coloring. We list the
colorings up to $r = 5$, \;\ as examples, below. The program is available
for download from the author's web site.\\

\noindent \textbf{Definitions:}\\
$C$ = list of the coloring in order. \\
$L$ = length of each interval
(proportional to each other) corresponding to each color in $C$. \\
$N$ = number of monochromatic Schur triple according to $C$ and $L$.

\[\begin{array}{llll} r = 1,& C = [1],& L = [1],& \mbox{N} =
\frac{n^2}{4} +O(n).  \\ r = 2,& C = [1,2,1],& L =
[1,\frac{3}{2},\frac{1}{4}],& \mbox{N} = \frac{n^2}{22} +O(n).\\ r =
3,& C = [1,2,1,3,1,2,1],& L = [1,\frac{3}{2},\frac{1}{4},3,\frac{1}{8},
\frac{487}{440}, \frac{47}{440}],& \mbox{N} = \frac{47n^2}{6238}
+O(n) \\ &&& \hspace{5.5mm} \sim \frac{n^2}{132.7234} +O(n).\\ \end{array}
\]
\noindent For $r \geq 4$, the lengths of the intervals are fractions with huge
numerators and denominators. So we omit $C$ and $L$ here.\\

\noindent $r =4$, N =
$\frac{69631222699293042329481527n^2}{67076984091396704809405315398}+O(n)
\sim \frac{n^2}{963.3176} +O(n) $.\\
$r =5$, N $\sim \frac{n^2}{7610.0730} +O(n) $.\\

\noindent For $r =6$, the lengths of the intervals are even larger
fractions. This caused Maple to slow down. We waited for about 8 hours
and we stopped. We did not get an answer.  However we were not really
disappointed about this failure. The algorithm is more important. \\

\noindent \textbf{4.2 Lower bounds} \\

The method used to find a lower bound in the previous two
sections could not be adapted for $r$-colorings, $r \geq 3$.
We did not make any progress for a lower bound of $r$-coloring cases.


\section{Conclusion}

We have new upper bounds for
triples $x+ay=z, a \geq 2$, in the $2$-coloring case. We also
have new upper bounds for Schur triples $x+y=z$, for $r$-colorings, $r \geq 3$
that considerably improve those of
\cite{LR}. But we failed to match the lower
and upper bounds for these two problems. There is
a possibility that other arguments in other papers \cite{BD}, \cite{RZ}
and \cite{Schoen} for the lower bound used in the original
problem can be adapted for the $r$-coloring problem. But the details of such an
argument seem complicated. We believe these upper bounds are actually
optimal. There might even be a beautiful simple way to solve it, but we
failed to find one (if it exists). We leave them as conjectures.\\

\noindent \textbf{Conjectures}: \\

\noindent 1) The (asymptotic) number of minimum monochromatic triples of
the form $\{x,y,x+ay\}, \;\ a \geq 2$ of 2-colorings of $[1,n]$,
are $\frac{n^2}{2a(a^2+2a+3)}+O(n)$. \\

\noindent 2) The (asymptotic) number of minimum Schur triples of
$r$-colorings of $[1,n]$, $r \geq 3$, are the same as the upper bounds
obtained from the Greedy Algorithm.

\subsection*{Acknowledgement}
I want to thank Bruce M. Landman and Aaron Robertson for their
beautiful book that revealed Ramsey theory to all of us. I really
enjoyed reading this book. I want to thank my advisor, Doron
Zeilberger, for teaching me to do symbolic programming and making the
graduate study years so much fun. I also want to thanks Yoni Berkowitz
for helpful discussions. Also thanks to my twin brother, Thotsaphon
Thanatipanonda, for helping with the programming at the beginning
phase of solving these problems.
\section*{Appendix}  
\appendix{}

\section{About the program}
\noindent $LowerBound(k,C)$ \\
\underline{input}: \;\ the number of intervals $k$, \;\
list of types of upper bound $C$ of $T_{i,k-i+1}$.\\
\noindent \underline{output}: lower bound of $M_{\chi}(n)$,
\;\ the upper bound of $Q$ and the optimal solution of $Q$.\\

\noindent $LowerBound2(k,C1,C2,a)$ \\
\underline{input}: \;\ the number of intervals $k$, \;\
list of types of upper bound $C1$ and $C2$ of $T_{i,k-i+1}$
and number $a$ in equation $x+ay=z$.\\
\noindent \underline{output}: lower bound of $M_{\chi,a}(n)$,
\;\ the upper bound of $Q_a$ and the optimal solution of $Q_a$.\\

\noindent $minAllST(n,r)$ \\
\underline{input}: \;\ length of intervals $n$, number of colors $r$.\\
\noindent \underline{output}: the r-coloring of all the interval of length n that
has the least number of monochromatic Schur triples.\\

\noindent $Ord(C,L,n)$\\
\underline{input}: \;\ the list of coloring, the
list of length corresponding to each color in $C$, symbol $n$.\\
\noindent \underline{output}: the number of the monochromatic Schur triples of order
$n^2$.\\

\noindent $Zeil(r)$\\
\noindent \underline{input}: \;\ number of color $r$.\\
\noindent \underline{output}: the coloring with length of
each coloring and also the total number of triples of order $n^2$
obtained from the Greedy Algorithm.

\end{document}